\documentclass[twoside]{article}

\usepackage[accepted]{aistats2020}
\usepackage{graphicx}
\usepackage{flushend}
\usepackage{latexsym}
\usepackage{cite}
\usepackage{amssymb}
\usepackage{float}
\usepackage{bbm}
\usepackage{textpos}
\usepackage{enumerate}
\usepackage{bm}
\usepackage{mathtools}
\usepackage[usenames,dvipsnames]{xcolor}
\usepackage{stmaryrd}
\usepackage{amsmath, amssymb,accents}
\usepackage{dsfont}
\usepackage{algorithmicx}
\usepackage{algorithm} 
\usepackage[noend]{algpseudocode} %
\usepackage{comment}
\usepackage{accents}
\usepackage{setspace}
\usepackage{theorem}
\usepackage{mathtools}
\usepackage{fancyhdr}
\allowdisplaybreaks

\newtheorem{definition}{Definition}

\newcommand{\supp}{\mathsf{supp}}
\newtheorem{theorem}{Theorem}
\newtheorem{lemma}{Lemma}

\newtheorem{corollary}{Corollary}

\usepackage{tikz}
\usetikzlibrary{calc}

\newcommand\numberthis{\addtocounter{equation}{1}\tag{\theequation}}
\newcommand{\sF}{\mathsf{F}}
\newcommand{\sS}{\mathsf{S}}

\newcommand{\cF}{\mathcal{F}}

\newcommand{\cK}{\mathcal{K}}

\newcommand{\cN}{\mathcal{N}}

\newcommand{\cP}{\mathcal{P}}

\newcommand{\cS}{\mathcal{S}}

\newcommand{\cX}{\mathcal{X}}

\newcommand{\EE}{\mathbb{E}}

\newcommand{\NN}{\mathbb{N}}

\newcommand{\PP}{\mathbb{P}}

\newcommand{\RR}{\mathbb{R}}

\newcommand*{\dd}{\, \mathsf{d}}
\newcommand{\vasti}{\bBigg@{3.5 }}
\newcommand{\vast}{\bBigg@{4}}
\newcommand{\Vast}{\bBigg@{5}}
\newcommand{\Vastt}{\bBigg@{7}}
\makeatother
\newcommand{\be}{\begin{equation}}
\newcommand{\ee}{\end{equation}}
\newcommand{\ba}{\begin{align}}
\newcommand{\ea}{\end{align}}
\newcommand{\baa}{\begin{align*}}
\newcommand{\eaa}{\end{align*}}

\newcommand*{\Noise}{\mathcal{G}_\sigma}
\newcommand*{\noise}{g_\sigma}
\newcommand*{\noised}{\tilde{g}_\sigma}
\newcommand*{\gauss}{\varphi_\sigma}
\newcommand*{\gaussd}{\tilde{\varphi}_\sigma}

\newcommand*{\Gaussk}{\mathcal{N}_{\sigma_k}}
\newcommand*{\Gauss}{\mathcal{N}_\sigma}
\newcommand*{\Gaussone}{\mathcal{N}_{\sigma_1}}
\newcommand*{\Gausstwo}{\mathcal{N}_{\sigma_2}}
\newcommand*{\Gaussmix}{\mathcal{N}_{\sqrt{\sigma_2^2-\sigma_1^2}}}

\newcommand{\var}{\mathsf{var}}

\newcommand{\wass}{\mathsf{W}_1}

\newcommand{\gwass}{\mathsf{W}_1^{(\sigma)}}
\newcommand{\Gwass}{\mathsf{W}_1^{(\mathcal{G}_\sigma)}}
\newcommand{\gwassone}{\mathsf{W}_1^{(\sigma_1)}}
\newcommand{\gwasstwo}{\mathsf{W}_1^{(\sigma_2)}}

\newcommand{\empmu}{\hat{\mu}_n}
\newcommand{\empnu}{\hat{\nu}_n}
\newcommand{\prodmu}{\mu^{\otimes n}}

\usepackage{caption}
\usepackage{subcaption}
\makeatletter

\begin{document}

%
\runningtitle{Gaussian-Smooth Optimal Transport: Metric Structure and Statistical Efficiency}

%

\twocolumn[

\aistatstitle{Gaussian-Smoothed Optimal Transport:\\
Metric Structure and Statistical Efficiency}

\aistatsauthor{Ziv Goldfeld \And Kristjan Greenewald}

\aistatsaddress{Cornell University \And IBM Research}
]

\begin{abstract}
Optimal transport (OT), and in particular the Wasserstein distance, has seen a surge of interest and applications in machine learning. However, empirical approximation under Wasserstein distances suffers from a severe curse of dimensionality, rendering them impractical in high dimensions. As a result, entropically regularized OT has become a popular workaround. However, while it enjoys fast algorithms and better statistical properties, it looses the metric structure that Wasserstein distances enjoy.  
This work proposes a novel Gaussian-smoothed OT (GOT) framework, that achieves the best of both worlds: preserving the 1-Wasserstein metric structure while alleviating the empirical approximation curse of dimensionality. Furthermore, as the Gaussian-smoothing parameter shrinks to zero, GOT $\Gamma$-converges towards classic OT (with convergence of optimizers), thus serving as a natural extension. An empirical study that supports the theoretical results is provided, promoting Gaussian-smoothed OT as a powerful alternative to entropic OT.
\end{abstract}

\section{Introduction}\label{SEC:intro}

In recent years optimal transport (OT) has been applied to a host of machine learning (ML) tasks as a powerful means of comparing probability measures. The Kantorovich OT \cite{kantorovich1942translocation} problem between two probability measures $\mu$ and $\nu$ with cost $c(x,y)$ is given by
\begin{equation}
    \inf_{\pi\in\Pi(\mu,\nu)}\int c(x,y)\dd\pi(x,y),\label{EQ:Kantorovich_OT}
\end{equation}
where $\Pi(\mu,\nu)$ is the set of transport plans (or couplings) between $\mu$ and $\nu$. Applications of the Kantorovich formulation include data clustering \cite{ho2017multilevel}, density ratio estimation \cite{iyer2014maximum}, domain adaptation \cite{courty2016optimal,courty2014domain}, generative models \cite{arjovsky2017wasserstein,gulrajani2017improved}, image recognition \cite{rubner2000earth,sandler2011nonnegative,li2013novel}, word and document embedding \cite{alvarez2018gromov,yurochkin2019hierarchical,grave2019unsupervised}, and many others. 

This surge in popularity has been driven by some highly advantageous properties of OT. Beyond its robustness to mismatched supports of $\mu$ and $\nu$ (crucial for learning generative models), when $c(x,y)=\|x-y\|$, \eqref{EQ:Kantorovich_OT} becomes the 1-Wasserstein distance\footnote{Any $p$-Wasserstein distance has these properties.}, which (i) has the operational interpretation of minimizing work (or expected cost); (ii) metrizes weak (also known as, weak*) convergence of probability measures; and (iii) defines a constant speed geodesic in the space of probability measures (giving rise to a natural interpolation between measures). These advantages, however, come with a price as OT is generally hard to compute and suffers from the so-called curse of dimensionality.

Specifically, suppose we have $n$ independent samples $(X_i)_{i=1}^n$ from a Borel probability measure $\mu$ on $\RR^d$. Consider the fundamental question of how quickly the empirical measure $\hat{\mu}_n\triangleq \frac{1}{n}\sum_{i=1}^n\delta_{X_i}$ approaches $\mu$ in the 1-Wasserstein distance, i.e., the $\EE\wass(\empmu,\mu)$ rate of decay. This quantity is at the heart of empirical approximation under $\wass$ since it controls the error in various additional approximation setups, such as  $\EE\big|\wass(\empmu,\nu)-\wass(\mu,\nu)\big|$ (one-sample goodness of fit test), $\EE\big|\wass(\empmu,\empnu)-\wass(\mu,\nu)\big|$ (two-samples tests)\footnote{Note that while Wasserstein-type GANs in practice typically use the two-sample setup since the generator distribution is intractable to compute, fundamentally the GAN actually corresponds to a one-sample setup since infinite samples can be obtained from the generator network.}, and others; see \cite{panaretos2019statistical} for a review on statistical applications of the Wasserstein distance.
Since $\wass$ metrizes weak convergence \cite[Cor. 6.18]{villani2008optimal}, the Glivenko-Cantelli theorem \cite{varadarajan1958convergence} implies $\wass(\empmu,\mu)\to 0$ as $n\to\infty$. Unfortunately, the convergence \emph{rate} in $n$ drastically deteriorates with dimension, scaling at best as $n^{-\frac{1}{d}}$ for any measure $\mu$ that is absolutely continuous with respect to (w.r.t.) the Lebesgue measure \cite{dudley1969speed}. Note that the $n^{-\frac{1}{d}}$ rate is sharp for all $d>2$ (see \cite{dobric1995asymptotics} for sharper results). This renders empirical approximation under the Wasserstein distance infeasible in high dimensions -- a disappointing shortcoming given the dimensionality of data in modern ML tasks.

In light of the above, entropic OT emerged as an appealing alternative to Kantorovich OT. Its popularity has been driven both by algorithmic advances \cite{cuturi2013sinkhorn,altschuler2017near} and some better statistical properties it possesses \cite{genevay2016stochastic,montavon2016wasserstein,rigollet2018entropic}. Entropic OT regularizes the expected cost by a Kullback-Leibler (KL) divergence, forming:
\begin{equation}
    \sS_c^{(\epsilon)}(\mu,\nu)\triangleq\inf_{\pi\in\Pi(\mu,\nu)}\int c(x,y)\dd\pi(x,y) + \epsilon \mathsf{D}(\pi\|\mu\times\nu),
\end{equation}
where $c(x,y)$ is the cost and $\mathsf{D}(\alpha\|\beta)\triangleq\int \log\left(\frac{\dd\alpha}{\dd\beta}\right)\dd\alpha$ if $\alpha\ll\beta$ and $+\infty$ otherwise. While the Wasserstein distance suffers from the curse of dimensionality, \cite{genevay2019sample} showed that if $c$ is Lipschitz and infinitely differentiable, then $\EE\big|\sS_c^{(\epsilon)}(\empmu,\empnu)-\sS_c^{(\epsilon)}(\mu,\nu)\big|\in O\left(n^{-\frac{1}{2}}\right)$, in all dimensions (see \cite{mena2019statistical} for sharper results 
specialized to quadratic cost). Despite this fast convergence in the two-sample test, sample complexity bounds in the (stronger) one-sample regime are not available. More importantly, the assumptions from \cite{genevay2019sample} exclude the distance cost $c(x,y)=\|x-y\|$, which is our main interest. 
Another drawback is that $\sS_c^{(\epsilon)}(\mu,\nu)$ is not a metric, even when $c(x,y)$ is \cite{feydy2018interpolating,bigot2019central} (e.g, $\sS_c^{(\epsilon)}(\mu,\mu)\neq 0$).\footnote{$\sS_c{(\epsilon)}$ can be transformed into a Sinkhorn divergence for which $\sS_c^{(\epsilon)}(\mu,\mu)=0)$, but it still is not a metric \cite{bigot2019central} since it lacks the triangle inequality.} Hence entropic OT retains several gaps in statistical convergence guarantees, and more importantly, it surrenders desirable properties of the Wasserstein distance. We thus seek an alternative OT framework that enjoys the best of both worlds.



\textbf{Contributions.} This paper proposes a novel OT framework, termed Gaussian-smoothed OT (GOT) 
that inherits the metric structure of $\wass$ while attaining stronger statistical guarantees than available for entropic OT. GOT of parameter $\sigma\geq 0$ between two $d$-dimensional probability measures $\mu$ and $\nu$ is defined~as
\begin{equation}
    \gwass(\mu,\nu)\triangleq \wass(\mu\ast\Gauss,\nu\ast\Gauss),\label{EQ:GOT}
\end{equation}
where $\ast$ stands for convolution and $\Gauss\triangleq\cN(0,\sigma^2\mathrm{I}_d)$ is the isotropic Gaussian measure of parameter $\sigma$. In other words, $\gwass(\mu,\nu)$ is simply the $\wass$ distance between $\mu$ and $\nu$ after each is smoothed by an isotropic Gaussian kernel.

We first show that just as $\wass$, for any fixed $\sigma\in[0,+\infty)$, $\gwass$ is a metric on the space of probability measures that metrizes the weak topology. Namely, a sequence of probability measures $(\mu_k)_{k\in\NN}$ converges weakly to $\mu$ if and only if $\gwass(\mu_k,\mu)\to 0$. We then turn to study properties of $\gwass(\mu,\nu)$ as a function of $\sigma$ for fixed $\mu$ and $\nu$. We establish continuity and non-increasing monotonicity. These, in particular, imply convergence of the optimal transportation costs, i.e., $\lim_{\sigma\to 0}\gwass(\mu,\nu)=\wass(\mu,\nu)$. Additionally, using the notion of $\Gamma$-convergence \cite{dalmaso2012Gamma_conv}, we establish convergence of optimizing transport plans. Thus, if $(\pi_{k})_{k\in\NN}$ is sequence of optimal transport plans for $\wass^{(\sigma_k)}(\mu,\nu)$, where $\sigma_k\to 0$, then $(\pi_{k})_{k\in\NN}$ converges weakly to an optimal plan for $\wass(\mu,\nu)$.

Lastly, we explore the one-sample empirical approximation under GOT, i.e., the convergence rate of $\EE\gwass(\empmu,\mu)$. It was shown in \cite{Goldfeld2019convergence} that Gaussian smoothing alleviates the curse of dimensionality, with $\EE\gwass(\empmu,\mu)$ converging as $n^{-\frac{1}{2}}$ in all dimensions. Although GOT is specialized to Gaussian noise, we present a generalized empirical approximation result that accounts for any subgaussian noise density. This, in turn, implies fast convergence of $\EE\big|\gwass(\empmu,\nu)-\gwass(\mu,\nu)\big|$ and  $\EE\big|\gwass(\empmu,\empnu)-\gwass(\mu,\nu)\big|$ via the triangle inequality. The expected value analysis is followed by a high probability claim derived through McDiarmid's inequality. Numerical results that validate these theoretical findings are provided. 
We conclude that GOT is an appealing alternative to entropic optimal transport, both in terms of its analytic and its statistical properties.

\section{Notation and Preliminaries}

Let $\cP(\RR^d)$ be the set of Borel probability measures on $\RR^d$, while $\cP_1(\RR^d)\subset\cP(\RR^d)$ are those with finite first moments, i.e.,  $\int_{\RR^d}\|x\|\dd\mu(x)<\infty$, where $\|\cdot\|$ is the Euclidean norm. We denote by $\Pi(\mu,\nu)\subset\cP(\RR^d)$ the set of transport plans (or couplings) between measures $\mu,\nu\in\cP(\RR^d)$. Namely, any $\pi\in\Pi(\mu,\nu)$ is a probability measure on $\RR^d\times\RR^d$ whose first and second marginals are $\mu$ and $\nu$, respectively.

The $n$-fold product extension of $\mu\in\cP(\RR^d)$ is $\prodmu$. 
The probability density function (PDF) of the isotropic Gaussian measure $\mathcal{N}_\sigma$ is $\gauss$. Given  $\mu,\nu\in\cP(\RR^d)$, their convolution $\mu\ast \nu\in\cP(\RR^d)$ is $(\mu\ast \nu)(\mathcal{A})=\int\int\mathds{1}_{\mathcal{A}}(x+y)\dd \mu(x)\dd \nu(y)$, where $\mathds{1}_\mathcal{A}$ is the indicator of $\mathcal{A}$. For two independent random variables $X\sim\mu$ and $Y\sim\nu$, we have $X+Y\sim\mu\ast\nu$.

We use $\EE_\mu f$ for the expectation of a measurable $f$ w.r.t. $\mu$, sometimes writing $\EE_\mu f(X)$ to emphasize its dependence on $X\sim\mu$. When the underlying probability measure is clear from the context, the subscript is omitted. Accordingly, the characteristic function of $\mu\in\cP(\RR^d)$ is $\phi_\mu(t)\triangleq\EE_{\mu}\big[e^{it^\top X}\big]$. For any $\mu,\nu\in\cP(\RR^d)$, we have $\phi_{\mu\ast\nu}(t)=\phi_\mu(t)\phi_\nu(t)$; if $\mu\times\nu\in\cP(\RR^d\times\RR^d)$ is the product measure of $\mu$ and $\nu$, then $\phi_{\mu\times\nu}(t,s)=\phi_\mu(t)\phi_\nu(s)$.


\begin{definition}[Weak Topology]
The \emph{weak topology} on $\cP(\RR^d)$ is induced by integration against the set $C_b^0(\RR^d)$ of bounded and continuous functions. Accordingly, we say that $(\mu_k)_{k\in\NN}\subset\cP(\RR^d)$ converges weakly to $\mu\in\cP(\RR^d)$, denoted by $\mu_k\rightharpoonup \mu$, if $\int_{\RR^d}f(x)\dd\mu_k(x)\to\int_{\RR^d}f(x)\dd\mu(x)$, for all $f\in C_b^0(\RR^d)$.
\end{definition}

It is a well-known fact that $\left(\cP_1(\RR^d),\wass\right)$ is a metric space, and that the 1-Wasserstein distance metrizes the weak topology (cf. \cite[Thm. 6.9]{villani2008optimal}). As shown in the sequel, this statement remains true if the 1-Wasserstein distance is replaced with its Gaussian-smoothed version, as defined next.

\begin{definition}[Gaussian-Smoothed $\bm{\wass}$] The Gaussian-smoothed 1-Wasserstein distance between $\mu,\nu\in\cP_1(\RR^d)$ is $\gwass(\mu,\nu)\triangleq \wass(\mu\ast\Gauss,\nu\ast\Gauss)$.
\end{definition}

Letting $X\sim\mu$, $Y\sim\nu$ and $Z,Z'\sim\Gauss$ be independent random variables, $\gwass(\mu,\nu)$ is the 1-Wasserstein distance between the probability laws of $X+Z\sim\mu\ast\Gauss$ and $Y+Z'\sim\nu\ast\Gauss$. Thus, $\gwass(\mu,\nu)$ can be understood as a `smoothed' version of $\wass$, where `smoothing' is applied to the probability measures via convolution with a Gaussian kernel (or, equivalently, via additive white Gaussian noise).

The theoretical results in this paper are organized as follows. Section \ref{SEC:metric_prop} studies the metric properties of $\gwass$. Section \ref{SEC:sigma_prop} establishes properties of $\gwass$ as a function of $\sigma$. One-sample empirical approximation rates under $\gwass$ are explored in Section \ref{SEC:stat_prop}.



\section{Metrizing the Weak Topology}\label{SEC:metric_prop}

Clearly, $\gwass(\mu,\nu)<+\infty$, for any $\mu,\nu\in\cP_1(\RR^d)$. Furthermore, similar to the regular 1-Wasserstein distance, $\gwass$ is a metric on $\cP_1(\RR^d)$, whose convergence is equivalent to convergence in the weak topology. 

\begin{theorem}[GOT Metric]\label{TM:Gwass_metric}
For any $\sigma\geq 0$, $\gwass:\cP_1(\RR^d)\times\cP_1(\RR^d)\to [0,+\infty)$ is a metric on $\cP_1(\RR^d)$.
\end{theorem}
This result mostly follows from $\wass$ being a metric. Some work is needed to establish the `identity of indiscernibles' properties. See Section \ref{SUBSEC:Gwass_metric_proof} for the proof.

\begin{theorem}[Weak Topology Metrization]\label{TM:Gwass_weak_top} Let $\sigma\geq 0$, $(\mu_k)\subset\cP(\RR^d)$ and $\mu\in\cP(\RR^d)$. Then $\gwass(\mu_k,\mu)\to 0$ if and only if (iff) $\EE_{\mu_k}\|X\|\to\EE_{\mu}\|X\|$ and $\mu_k\rightharpoonup \mu$. Consequently, $\gwass(\mu_k,\mu)\to 0$ iff $\wass(\mu_k,\mu)\to 0$.
\end{theorem}

Theorem \ref{TM:Gwass_weak_top} with $\wass$ in place of $\gwass$ is a well-known result \cite[Thm. 6.9]{villani2008optimal}). The above can be therefore understood as the statement that `the 1-Wasserstein topology is invariant to convolutions with Gaussian kernels'. See Section \ref{SUBSEC:Gwass_weak_top_proof} for the proof.


\section{Dependence on Noise Parameter}\label{SEC:sigma_prop}

We study properties on $\gwass(\mu,\nu)$, for fixed $\mu,\nu\in\cP_1(\RR^d)$, as a function of $\sigma\in[0,+\infty)$. 

\begin{theorem}[GOT Dependence on $\bm{\sigma}$]\label{TM:sigma_prop}
Fix $\mu,\nu\in\cP_1(\RR^d)$. The following hold:
\begin{enumerate}[i)]
    \item $\gwass(\mu,\nu)$ is continuous and monotonically non-increasing in $\sigma\in[0,+\infty)$;
    \item $\lim_{\sigma\to 0}\gwass(\mu,\nu)=\wass(\mu,\nu)$;
    \item $\lim_{\sigma\to \infty}\gwass(\mu,\nu)\neq 0$, for some $\mu,\nu\in\cP_1(\RR^d)$.
\end{enumerate}
\end{theorem}
While $\gwass(\mu,\nu)$ is a monotonically non-increasing function of $\sigma$, as $\sigma \rightarrow \infty$ it is interestingly not true in general that $\wass(\mu\ast\Gauss,\nu\ast\Gauss)$ decays to zero. The proof of Theorem \ref{TM:sigma_prop} (Section \ref{SUBSEC:sigma_prop_proof}) shows this via a simple Dirac measure example.

A key technical tool (that may be of independent interest) for establishing item (i) above is the following lemma, which ties GOT at different noise levels to one another. Its proof (Section \ref{SUBSEC:Noise_std_relation_proof}) uses the Kantorovich-Rubinstein duality.

\begin{lemma}[Stability Across $\sigma$]\label{LEMMA:Noise_std_relation}
Fix $\mu,\nu\in\cP_1(\RR^d)$, and $0\leq\sigma_1<\sigma_2<+\infty$. We have
\begin{equation*}
    \gwasstwo(\mu,\nu)\mspace{-2mu}\leq\mspace{-2mu}\gwassone(\mu,\nu)\mspace{-2mu}\leq\mspace{-2mu}\gwasstwo(\mu,\nu)\mspace{-2mu}+2\sqrt{\mspace{-2mu}d\left(\sigma_2^2\mspace{-2mu}-\mspace{-2mu}\sigma_1^2\right)}.
\end{equation*}
\end{lemma}

Theorem \ref{TM:sigma_prop} established convergence of transport costs, i.e., that $\wass^{(\sigma_k)}(\mu,\nu)\to\gwass(\mu,\nu)$ as $\sigma_k\to\sigma$. The next result shows we also have convergence of optimal plans. Namely, a sequence of optimal couplings $(\pi_k)_{k\in\NN}$ for $\wass^{(\sigma_k)}(\mu,\nu)$ (weakly) approaches an optimal coupling for $\gwass(\mu,\nu)$ as $k$ goes to infinity.

\begin{theorem}[Convergence of Optimal Plans]\label{TM:convergence_couplings} Fix $\mu,\nu\in\cP_1(\RR^d)$ and let $(\sigma_k)_{k\in\NN}$ be a sequence~with $\sigma_k\searrow\sigma\geq 0$. 
Let $\pi_k\in\Pi(\mu\ast\Gaussk,\nu\ast\Gaussk)$, $k\in\NN$, be an optimal coupling for $\wass^{(\sigma_k)}(\mu,\nu)$. Then there exists  $\pi\in\Pi(\mu\ast\Gauss,\nu\ast\Gauss)$ such that $\pi_k\rightharpoonup \pi$ (weakly) as $k\to\infty$ and 
$\pi$ is optimal for $\gwass(\mu,\nu)$.
\end{theorem}

The proof of Theorem \ref{TM:convergence_couplings} (Section \ref{SUBSEC:convergence_couplings_proof}) relies on the notion of $\Gamma$-convergence. Convergence of optimal transport plans then follows by standard tightness arguments. In particular, this theorem implies that a sequence of optimal transport plans for $\gwass(\mu,\nu)$ converges to an optimal plan for the regular 1-Wasserstein distance $\wass(\mu,\nu)$ as $\sigma \rightarrow 0$.


\section{Empirical Approximation}\label{SEC:stat_prop}

We now explore statistical properties of $\gwass$. In fact, our derivation accounts for any isotropic noise distribution $\Noise$ that along each coordinate is $\sigma$-subgaussian with a bounded and monotone (in a proper sense) density.\footnote{A further extension to nonisotropic noise is possible via similar techniques, but we do not delve into it here.} Gaussian noise is captured as a special case. 

Consider the fundamental one-sample empirical approximation, where $\mu\in\cP_1(\RR^d)$ is approximated by $\empmu\triangleq\frac{1}{n}\sum_{i=1}^n\delta_{X_i}$, with $(X_1,\ldots,X_n)\sim \prodmu$ and $\delta_x$ as the Dirac measure centered at $x$. We study how fast $\Gwass(\empmu,\mu)\triangleq\wass(\empmu\ast\Noise,\mu\ast\Noise)\to 0$ with~$n$.\footnote{Of course, $\mathsf{W}_1^{(\Gauss)}(\mu,\nu)=\gwass(\mu,\nu)$.} In a remarkable contrast to the 1-Wasserstein curse of dimensionality, we show $\EE_{\prodmu}\gwass(\empmu,\mu)\in O\big(n^{-\frac{1}{2}}\big)$ in all dimensions \cite{Goldfeld2019convergence}, thus attaining the parametric rate.



To state the results, we first define subgaussianity.

\begin{definition}[Subgaussian Measure]\label{DEF:SG}
A probability measure $\mu\in\cP_1(\RR^d)$ is $K$-subgaussian, for $K>0$, if for any $\alpha \in \RR^d$, $X\sim \mu$ satisfies
\begin{equation}
\EE_\mu\Big[e^{\alpha^T(X -\EE X)}\Big] \leq e^{\frac{1}{2}K^2 \|\alpha\|^2}.\label{EQ:SG}
\end{equation}
\end{definition}

We begin with a bound on the expected value and then move to a high probability bound. The next theorem generalizes \cite[Prop. 1]{Goldfeld2019convergence} to non-Gaussian noise models.



\begin{theorem}[GOT Empirical Approximation]\label{TM:gwass_emp_expectation}
Fix $d\geq 1$, $\sigma>0$ and $K>0$. Let $\Noise\in\cP_1(\RR^d)$ have a density $\noise$ that decomposes as $\noise(x) = \prod_{j= 1}^d \noised(x_j)$. Assume that $\noised$ is $\sigma$-subgaussian, bounded and monotonically decreases as its argument goes away from zero in either direction. For any $K$-subgaussian $\mu\in\cP_1(\RR^d)$, we have
\begin{equation}
\EE_{\mu^{\otimes n}} \Gwass(\empmu,\mu) \leq c_{\sigma,d,K} n^{-\frac{1}{2}},\label{EQ:gwass_emp_expectation_equation}
\end{equation}
where $c_{\sigma,d,K}=e^{O(d)}$ is given in \eqref{EQ:W1_bound_constant}. In particular $\gwass(\empmu,\mu)\in O\left(n^{-\frac{1}{2}}\right)$.
\end{theorem}
The proof of Theorem \ref{TM:gwass_emp_expectation} is given in Section \ref{SUBSEC:gwass_emp_expectation_proof}.


\begin{corollary}[Concentration Inequality]\label{COR:gwass_emp_whp}
Under the paradigm of Theorem \ref{TM:gwass_emp_expectation}, denote $\cX\triangleq\supp(\mu)$ and suppose $\mathsf{diam}(\cX)<\infty$, where $\mathsf{diam}(\cX)=\sup_{x\neq y\in\cX}\|x-y\|$. For any $t>0$ we have
\begin{equation}
\PP_{\mspace{-3mu}\mu^{\mspace{-2mu}\otimes\mspace{-1mu} n}}\mspace{-4mu}\left(\Big|\Gwass\mspace{-2mu}(\empmu,\mu)\mspace{-1.5mu}-\mspace{-1.5mu}\EE\Gwass\mspace{-2mu}(\empmu,\mu)\Big|\mspace{-2mu}\geq\mspace{-2mu} t\mspace{-1.5mu}\right)\mspace{-3mu}\leq 2e^{-\frac{2t^2n}{\mathsf{diam}(\cX)^2}}\label{EQ:gwass_emp_whp_equation}
\end{equation}
and consequently, 
\begin{equation}
\PP_{\mu^{\otimes n}}\left(\Gwass(\empmu,\mu)\in \omega\left(\frac{\log n}{\sqrt{n}}\right)\right)\leq \frac{1}{\mathsf{poly}(n)}.\label{EQ:gwass_emp_whp_equation1}
\end{equation}
\end{corollary}

The proof Theorem \ref{COR:gwass_emp_whp} is given in Section \ref{SUBSEC:gwass_emp_whp_proof}. It uses the $\wass$ duality and McDiarmid's inequality.


\section{Empirical Results}\label{SEC:empirical_results}


\begin{figure*}[htbp]
\centering
 \begin{subfigure}[t]{0.49\textwidth}
  \includegraphics[width=0.9\linewidth]{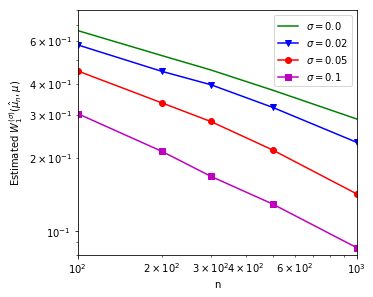}
   \caption{$d=5$}
 \end{subfigure} 
 \begin{subfigure}[t]{0.49\textwidth}
   \includegraphics[width=0.9\linewidth]{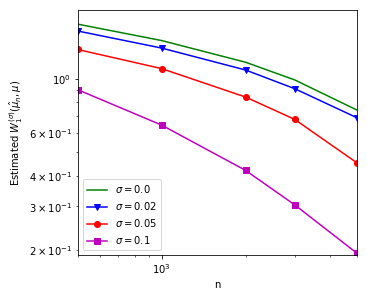}
   \caption{$d=10$}
 \end{subfigure}
  \caption{Convergence of $\gwass(\empmu,\mu)$ as a function of the number of samples $n$ for various values of $\sigma$, shown in log-log space. The measure $\mu$ is the uniform distribution over $[0,1]^d$. Note that $\sigma=0$ corresponds to the vanilla Wasserstein distance, which converges slower than GOT (note the difference in slopes), especially with larger $d$.}\label{FIG:convergence}
\end{figure*}


We turn to some numerical experiments demonstrating the difference in empirical approximation convergence rates between the regular 1-Wasserstein distance and GOT. Specifically, we compute $\wass(\empmu,\mu)$ and $\gwass(\empmu,\mu)$, for $\mu=\mathsf{Unif}\big([0,1]^d\big)$ the uniform measure on $[0,1]^d$, and $\empmu=\frac{1}{n}\sum_{i=1}^n\delta_{X_i}$ the empirical measure based on i.i.d. samples $X_1,\ldots,X_n\sim \mu$. This simple setup also hints at how broad the class of distributions for which $\wass(\empmu,\mu)$ attains the poor convergence rate. 

The GOT framework corresponds to the 1-Wasserstein distance between two continuous (smooth) distributions. To evaluate this 1-Wasserstein distance we chose to employ the neural network (NN) based dual optimization approach of \cite{gulrajani2017improved}. This approach seems to be better suited for continuous probability measures than, e.g., the Sinkhorm algorithm \cite{cuturi2013sinkhorn}. Starting from Kantorovich-Rubinstein duality

\begin{equation}
    \wass(\mu,\nu)=\sup_{\|f\|_{\mathsf{Lip}}\leq1}\EE_{\mu} f -\EE_{\nu} f,
\end{equation}
the function $f$ is first parametrized by a NN $f_\theta$, with parameter set $\theta\in\Theta$,\footnote{We used a fully connected network with 3 hidden ReLU layers, each comprising 1024 nodes. The network was trained until convergence of the estimated Wasserstein distance.} and then the $\|f_\theta\|_{\mathsf{Lip}\leq 1}$ constraint is relaxed to a regularization penalty on the expected gradient of $f_\theta(x)$ (w.r.t. to $x$). In sum, as in \cite{gulrajani2017improved}, we use the ADAM stochastic gradient ascent method to optimize
\begin{equation}
    \sup_{\theta\in\Theta} \EE_{\mu}f_\theta -\EE_{\nu}f_\theta +\lambda\EE_\eta \Big[\big(\|\nabla_x f_{\theta}\|-1\big)^2\Big],\label{EQ:Wass_GP}
\end{equation}
where $\eta$ interpolates between $\mu$ and $\nu$ in a manner compatible with the gradient penalty (GP) theoretical justification \cite[Prop. 1]{gulrajani2017improved}. Specializing to $\gwass$, $\mu$ and $\nu$ above are replaced with their Gaussian-smoothed versions, i.e., $\mu\ast\Gauss$ and $\nu\ast\Gauss$, respectively. To approximate expectations with empirical sums, we sample from these Gaussian-smoothed measures by adding (sampled) Gaussian noise to the original samples. This makes use of the fact that convolution of probability measures corresponds to sums of independent random variables.

Figure \ref{FIG:convergence} shows the results for $d = 5$ and $d = 10$, with each curve being the average of 10 random trials\footnote{Error bars were omitted since they were too small to be visible.}. Note the slower convergence to zero of the $\sigma=0$ case corresponding to vanilla $\wass$, compared to the approximately $O(n^{-1/2})$ convergence of the $W_1^{(\sigma)}$ metrics for larger $\sigma > 0$. In the $d=10$ plot, the curves slightly accelerate as $n$ increases instead of staying linear. This seems to originate from a two-fold imperfection in the NN-based approximation of the Lipschitz function $f$. First, the GP regularization does not perfectly enforce the Lipschitz constraint especially in high dimensions. Second, to accurately evaluate $\gwass(\empmu,\mu)$ the network effectively needs to overfit $\empmu$. As NNs tends to avoid overfitting (especially once the number of modes $n$ in $\empmu$ becomes large), additional slackness might be introduced.


As expected, the $\wass(\empmu,\mu)$ estimate converges significantly slower than its Gaussian-smoothed counterpart, as evident by comparing the slopes of the curves in log-log space. In particular, the convergence of the $\wass(\empmu,\mu)$ estimate is much slower for $d=10$ than for $d=5$ as predicted. The $\gwass$ estimate, on the other hand, still converges approximately as $O(n^{-1/2})$ in both cases. The fact that $\gwass$ is monotonically decreasing in $\sigma$ can also be seen from the plots. These results are comparable with the ones from \cite{genevay2019sample} for two-sample empirical approximation of entropic OT.

\section{Proofs}\label{SEC:Proofs}


\subsection{Proof of Theorem \ref{TM:Gwass_metric}}\label{SUBSEC:Gwass_metric_proof}

The fact that $\gwass(\mu,\nu)$ is symmetric, non-negative and equals zero when $\mu=\nu$ follows from its definition.

To prove the triangle inequality, i.e., $\gwass(\mu_1,\mu_3)\leq\gwass(\mu_1,\mu_2)+\gwass(\mu_2,\mu_3)$, for any $\mu_1,\mu_2,\mu_3\in\cP_1(\RR^d)$, let $\pi_{12}\in\Pi(\mu_1\ast\Gauss,\mu_2\ast\Gauss)$ and $\pi_{23}\in\Pi(\mu_2\ast\Gauss,\mu_3\ast\Gauss)$ be optimal couplings for $\gwass(\mu_1,\mu_2)$ and $\gwass(\mu_2,\mu_3)$, respectively. Applying the Gluing Lemma \cite{villani2008optimal}, let $\pi\in\cP_1(\RR^d\times\RR^d\times\RR^d)$ be a probability measure with $\pi_{12}$ and $\pi_{23}$ as its marginals on the corresponding coordinates. Letting $\pi_{13}(A\times B)\triangleq\pi(A\times\RR^d\times B)$, we have $\pi_{13}\in\Pi(\mu_1,\mu_3)$ and 
\begin{align*}
    \gwass(\mu_1,\mu_3)&\leq \mathbb{E}_{\pi_{13}}\|X_1-X_3\|\\
    &\leq \mathbb{E}_{\pi_{12}}\|X_1-X_2\|+\mathbb{E}_{\pi_{23}}\|X_2-X_3\|\\
    &=\gwass(\mu_1,\mu_2)+\gwass(\mu_2,\mu_3).\numberthis
\end{align*}

It remains to show that $\gwass(\mu,\nu)=0$ implies that $\mu=\nu$. Since $\wass$ is a metric, we know that if $\gwass(\mu,\nu)=0$ then $\mu\ast\Gauss=\nu\ast\Gauss$. This implies pointwise equality between characteristic functions: $\phi_\mu\phi_{\Gauss}=\phi_\nu\phi_{\Gauss}$. Since $\phi_{\Gauss}\neq 0$ everywhere, we get $\phi_\mu=\phi_\nu$ pointwise, implying $\mu=\nu$.


\subsection{Proof of Theorem \ref{TM:Gwass_weak_top}}\label{SUBSEC:Gwass_weak_top_proof}

The claim relies on the equivalence between weak convergence and pointwise convergence of characteristic functions. Since $\wass$ metrizes weak convergence:
\begin{align*}
    &\gwass(\mu_k,\mu)\to 0\\
    &\ \ \iff \quad \mu_k\ast\Gauss \to \mu\ast\Gauss\\
    &\ \ \iff \quad \phi_{\mu_k}(t)\phi_{\Gauss}(t)=\phi_\mu(t)\phi_{\Gauss}(t)\,,\quad \forall t\in\RR^d\\
    &\ \ \iff \quad \phi_{\mu_k}(t)=\phi_\mu(t)\,,\quad \forall t\in\RR^d.
\end{align*}

\subsection{Proof of Theorem \ref{TM:sigma_prop}}\label{SUBSEC:sigma_prop_proof}

For Claim (ii), the fact that $\lim_{\sigma\to0}\gwass(\mu,\nu)=\wass(\mu,\nu)$ follows from Lemma \ref{LEMMA:Noise_std_relation} by taking $\sigma_1=0$ and $\sigma_2=\sigma\to 0$.

For Claim (i), $\gwass(\mu,\nu)$ being monotonically non-increasing in $\sigma$ also follows directly from Lemma \ref{LEMMA:Noise_std_relation}. To prove continuity at $\sigma\in(0,+\infty)$, we consider left- and right- continuity separately. Let $\sigma_k\nearrow \sigma$ as $k\to\infty$. Lemma \ref{LEMMA:Noise_std_relation} gives
\begin{equation}
    \gwass(\mu,\nu)\leq\wass^{(\sigma_k)}(\mu,\nu)\leq\gwass(\mu,\nu)+2d\sqrt{\sigma^2-\sigma_k^2},
\end{equation}
and left-continuity follows.

To see that $\gwass(\mu,\nu)$ is right-continuous in $\sigma$, let $\sigma_k\searrow\sigma$ and denote $\epsilon_k\triangleq\sqrt{\sigma_k^2-\sigma^2}$. We have
\begin{align*}
    \wass^{(\sigma_k)}(\mu,\nu)=\wass^{(\epsilon_k)}(\mu\ast\Gauss,\nu\ast\Gauss)\xrightarrow[k\to\infty]{}\gwass(\mu,\nu),\numberthis
\end{align*}
where the last step uses $\gwass$ continuity at $\sigma=0$.

Moving to Claim (iii), let $\mu=\delta_x$ and $\nu=\delta_y$ be two Dirac measures at $x\neq y\in\RR^d$. For any $\sigma\in[0,+\infty)$, we have
\begin{align*}
    \gwass(\mu,\nu)&=\wass\big(\mathcal{N}(x,\sigma^2\mathrm{I}_d),\mathcal{N}(y,\sigma^2\mathrm{I}_d)\big)\\
    &\geq \Big\|\EE_{\mathcal{N}(x,\sigma^2\mathrm{I}_d)}X-\EE_{\mathcal{N}(y,\sigma^2\mathrm{I}_d)}Y\Big\|\\
    &=\|x-y\|,
\end{align*}
where the equality uses Jensen's inequality and convexity of norm.

\subsection{Proof of Lemma \ref{LEMMA:Noise_std_relation}}\label{SUBSEC:Noise_std_relation_proof}

The first inequality immediately follows because $\wass$ is non-increasing under convolutions and since $\Gausstwo=\Gaussone\ast\Gaussmix$. 

For the second inequality, we use Kantorovich-Rubinstein duality to write
\begin{align*}
    &\gwassone(\mu,\nu)=\sup_{\|f_1\|_{\mathsf{Lip}}\leq 1} \EE_{\mu\ast\Gaussone}f_1-\EE_{\nu\ast\Gaussone}f_1;\\
    &\gwasstwo(\mu,\nu)=\sup_{\|f_2\|_{\mathsf{Lip}}\leq 1} \EE_{\mu\ast\Gausstwo}f_2-\EE_{\nu\ast\Gausstwo}f_2.\\
\end{align*}
Letting $f_1^\star$ be optimal for $\gwassone(\mu,\nu)$, we have
\begin{equation}
    \gwasstwo(\mu,\nu)\geq \EE_{\mu\ast\Gausstwo}f_1^\star-\EE_{\nu\ast\Gausstwo}f_1^\star.\label{EQ:dual_lb}
\end{equation}
Set $X\sim\mu$, $Z_1\sim\Gaussone$ and $Z_{21}\sim\Gaussmix$ as independent random variables; clearly, $Z_2\triangleq Z_1+Z_{21}\sim\Gausstwo$. Consider:
\begin{subequations}
\begin{align*}
    \left|\EE_{\mu\ast\Gaussone}f_1^\star-\EE_{\mu\ast\Gausstwo}f_1^\star\right|&=\EE f_1^\star(X+Z_1)-\EE f_1^\star(X+Z_2)\\
    &\leq \EE\|Z_{21}\|\\
    &=\sqrt{d\left(\sigma_2^2-\sigma_1^2\right)},\numberthis\label{EQ:dual_lb1}
\end{align*}
where the last in equality uses $\|f_1^\star\|_{\mathsf{Lip}}\leq 1$. Similarly, one has
\begin{equation}
    \left|\EE_{\nu\ast\Gaussone}f_1^\star-\EE_{\nu\ast\Gausstwo}f_1^\star\right|\leq \sqrt{d\left(\sigma_2^2-\sigma_1^2\right)}.\label{EQ:dual_lb2}
\end{equation}\label{EQ:dual_lb12}%
\end{subequations}
Inserting \eqref{EQ:dual_lb12} into \eqref{EQ:dual_lb} concludes the proof.


\subsection{Proof of Theorem \ref{TM:convergence_couplings}}\label{SUBSEC:convergence_couplings_proof}

We first include the definitions of tightness of measures and $\Gamma$-convergence of functionals.

\begin{definition}[Tightness of Measures]
A subset $\cS\subset\cP(\RR^d)$ is tight if for any $\epsilon>0$ there is a compact set $\cK_\epsilon\subset\RR^d$ such that $\mu(\cK_\epsilon)\geq 1-\epsilon$, for all $\mu\in\cP(\RR)^d$.
\end{definition}

\begin{definition}[$\bm{\Gamma}$-Convergence]
Let $\cX$ be a metric space and $\sF_k:\cX\to\RR$, $k\in\NN$ be a sequence of functionals. We say $(\sF_k)_{k\in\NN}$ $\Gamma$-converges to $\sF:\cX\to\RR$, and we write $\sF_k\overset{\Gamma}{\to}\sF$, if:
\begin{enumerate}[i)]
    \item For every $x_k,x\in\cX$, $k\in\NN$, with $x_k\to x$, we have $\sF(x)\leq\liminf_{k\to\infty}\sF_k(x_k)$;
    \item For any $x\in\cX$, there exists $x_k\in\cX$, $k\in\NN$, with $x_k\to x$, and $\sF(x)\geq\limsup_{k\to\infty}\sF_k(x_k)$
\end{enumerate}
\end{definition}

By pointwise convergence of characteristic functions, $P_k\triangleq\mu\ast\Gaussk$ and $Q_k\triangleq\nu\ast\Gaussk$ are weakly convergent measures on $\RR^d$. Prokhorov's Theorem then implies they are tight. By \cite[Lemma 4.4]{villani2008optimal} we have that $\Pi\Big((P_k)_{k\in\NN},(Q_k)_{k\in\NN}\Big)$, the set of all couplings with marginals in $(P_k)_{k\in\NN}$ and $(Q_k)_{k\in\NN}$, is also tight. Hence, the sequence of optimal couplings $(\pi_k)_{k\in\NN}$ is tight and weakly converges to some $\pi\in\cP(\RR^d\times\RR^d)$. Taking the limit of the relation $\pi_k\in\Pi(P_k,Q_k)$ we obtain $\pi\in\Pi(P,Q)$, where $P\triangleq\mu\ast\Gauss$ and $Q\triangleq\nu\ast\Gauss$.

With that in mind, recall that if $(\sF_k)_{k\in\NN}$ $\Gamma$-converges to $\sF$, then $\lim_{k\to\infty}\inf\sF_k=\inf\sF$ \cite[Thm. 7.8]{dalmaso2012Gamma_conv}. Furthermore, if $(x_k)_{k\in\NN}$ is a sequence of minimizers of $\sF_k$, for each $k\in\NN$, then any cluster (limit) point of $(x_k)_{k\in\NN}$ is a minimizer of $\cF$ \cite[Cor. 7.20]{dalmaso2012Gamma_conv}. Thus, to conclude the proof of Theorem \ref{TM:convergence_couplings} it suffices to establish $\Gamma$-convergence of $\sF_k:\cP(\RR^d\times\RR^d)\to\RR\cup\{\infty\}$ to $\sF:\cP(\RR^d\times\RR^d)\to\RR\cup\{\infty\}$ defined as
\begin{align*}
    \sF_k(\pi)&=\begin{cases}\EE_\pi\|X-Y\|,\quad \pi\in\Pi(\mu\ast\Gaussk,\nu\ast\Gaussk)\\
    \infty,\quad\quad\quad\quad\quad \mbox{otherwise}\end{cases}\\
    \sF(\pi)&=\begin{cases}\EE_\pi\|X-Y\|,\quad \pi\in\Pi(\mu\ast\Gauss,\nu\ast\Gauss)\\
    \infty,\quad\quad\quad\quad\quad \mbox{otherwise}\end{cases}.\numberthis
\end{align*}

We start with the $\liminf$ $\Gamma$-convergence inequality. First observe that if $(\pi_k)_{k\in\NN}$ does not contain a subsequence (without relabeling) such that $\pi_k\in\Pi(\mu\ast\Gaussk,\nu\ast\Gaussk)$, then the claim is trivial. Accordingly, assume (again, up to extraction of subsequences) that $\pi_k\in\Pi(\mu\ast\Gaussk,\nu\ast\Gaussk)$, for all $k\in\NN$. Since $x\mapsto\|x\|$ is a non-negative and continuous, the $\liminf$ condition directly follows from the Portmanteau Theorem:
\begin{align*}
    \sF(\pi)&=\int\|x-y\|\dd\pi\\
    &\leq \liminf_{k\to\infty}\int \|x-y\|d\pi_k\\
    &=\liminf_{k\to\infty}\sF_k(\pi_k).\numberthis
\end{align*}

For the $\limsup$ let $\pi\in\Pi(\mu\ast\Gauss,\nu\ast\Gauss)$. For convenience, we use random variable notation. There exists a tuple $(X,Y,Z',Z'')$ with marginal distributions $X\sim\mu$, $Y\sim\nu$ and $Z',Z''\sim\Gauss$, such that $(X,Z')$ are independent, $(Y,Z'')$ are independent, and $(X+Z',Y+Z'')\sim\pi$. 

To construct the sequence $(\pi_k)_{k\in\NN}$, let $Z_k\sim\cN_{\sqrt{\sigma_k^2-\sigma^2}}$ be independent of $(X,Y,Z',Z'')$. Setting $\pi_k$ as the joint probability law of $(X+Z'+Z_k,Y+Z''+Z_k)$, we have $\pi_k\in\Pi(\mu\ast\Gaussk,\nu\ast\Gaussk)$, $k\in\NN$. Evaluating $\sF_k$ we obtain
\begin{align*}
    \sF_k(\pi_k)&=\EE\|X+Z'+Z_k-Y-Z''-Z_k\|\\
    &=\EE\|X+Z'-Y-Z''\|\\
    &=\sF(\pi),
\end{align*}
which in particular implies the $\limsup$ condition.


\subsection{Proof of Theorem \ref{TM:gwass_emp_expectation}}\label{SUBSEC:gwass_emp_expectation_proof}

The 1-Wasserstein distance is upper bounded by weighted total variation (TV) as follows \cite[Theorem~6.15]{villani2008optimal}:
\begin{equation}
\wass(\empmu\ast\Noise, \mu\ast\Noise) \leq \int_{\RR^d} \|t\| \big|r_n(t) - q(t)\big| \dd t,
\end{equation}
where $r_n$ and $q$ are the densities of $\empmu\ast\Noise$ and $\mu\ast\Noise$, respectively. The inequality is proved using the maximal TV coupling of $\empmu\ast\Noise$ with $\mu\ast\Noise$. 

Let $a > 0$ (to be specified later) and set $f_a:\RR^d \to \RR$ as the density of $\mathcal{N}\left(0,\frac{1}{2a} \mathrm{I}_d\right)$. By Cauchy-Schwarz, we have
\begin{align*}
&\EE_{\prodmu}\int_{\RR^d} \|t\| \big|r_n(t) - q(t)\big| \dd t\\
&\leq\mspace{-3mu} \left(\int_{\RR^d} \mspace{-4mu}\|t\|^2f_a(t)\dd t\right)^{\mspace{-6mu}\frac{1}{2}}\mspace{-5mu}\left(\mathbb{E}_{\prodmu}\mspace{-5mu}\int_{\RR^d}\mspace{-6mu} \frac{\big(q(t)\mspace{-1.5mu}-\mspace{-1.5mu}r_n(t)\big)^2}{f_a(t)} \dd t\right)^{\mspace{-4mu}\frac{1}{2}}\mspace{-7mu}.\numberthis\label{EQ:W1_CS_bound}
\end{align*}
The first term 
equals $\frac{d}{2a}$. Turning to the second integral, note that $r_n(t)=\frac{1}{n}\sum_{i=1}^n \noise(t-X_i)$, where $\{X_i\}_{i=1}^n$ are i.i.d. and $\EE_\mu \noise(t-X_i)=q(t)$. Using the definition of subgaussianity (Definition \ref{DEF:SG}), we have the following lemma (proven in Appendix \ref{APPEN:density_bound_proof}) that bounds $\noise$ everywhere in terms of the Gaussian density~$\gauss$.
\begin{lemma}\label{LEMMA:density_bound}
Let $\delta\triangleq\min\left\{1,\frac{1}{4\sigma^2}\right\}$. There exists a constant $c_1>0$ 
such that
\begin{equation}
    \noise(t)\leq c_1^d e^{\delta \|t\|^2}\gauss(t),\quad\forall t\in\RR^d.
\end{equation}
\end{lemma}

We now can bound the second integrand of \eqref{EQ:W1_CS_bound}:
\begin{align*}
    \EE_{\prodmu}\big(q(z)-r_n(z)\big)^2&=\var_{\prodmu}\big(r_n(z)\big)\\
    &= \var_{\prodmu}\left( \frac{1}{n}\sum_{i=1}^n \noise(z-X_i)\right)\\
    &= \frac{1}{n} \var_{\mu}\big(\noise(z-X)\big)\\
    &\leq \EE_{\mu}\noise^2(z-X)\\
    &\leq c_1^d{\delta}^{2d}\EE_{\mu} e^{2 \delta \|z-X\|^2 }\gauss^2(z-X)\\
    &\leq\frac{c_2^2}{n}\EE_\mu e^{ -\frac{1}{2\sigma^2} \|z - X\|^2},\numberthis
\end{align*}
with $c_2 = c_1^d (2 \pi \sigma^2)^{-d/2}$. This further implies
\begin{equation}
\int_{\RR^d}\EE_{\prodmu} \frac{\big(q(t)-r_n(t)\big)^2}{f_a(t)} \dd z \leq \frac{c_2}{n2^{d/2}} \EE \frac{1}{f_a(X+ Z  )},\label{EQ:W1_1/f_bound}
\end{equation}
where $X\sim \mu$ and $Z\sim\Gauss$ are independent. 

Starting from \eqref{EQ:W1_1/f_bound}, we finish the proof via steps similar to \cite{Goldfeld2019convergence}. Specifically, for $c_3\triangleq\left(\frac{\pi}{a}\right)^{\frac{d}{2}}$, it holds that $\big(f_a(t)\big)^{-1}= c_3 e^{a \|t\|^2}$. Since $X$ is $K$-subgaussian and $Z$ is $\sigma$-subgaussian, $X + Z$ is $(K+ \sigma)$-subgaussian. Following \eqref{EQ:W1_1/f_bound}, for any $0<a<\frac{1}{2(K + \sigma)^2}$, we have~\cite[Rmk. 2.3]{hsu2012tail}
\begin{align*}
&\frac{c_2}{n2^{d/2}}\EE \frac{1}{f(X + Z)}\\
&=\frac{c_2 c_3}{n2^{d/2}} \EE \exp\Big(a \big\|X + Z\big\|^2\Big)\numberthis\label{EQ:W1_int2_bound}\\
&\leq \frac{c_2 c_3}{n2^{d/2}} \exp\mspace{-3mu}\left(\mspace{-3mu}\big(K\mspace{-2mu}+\mspace{-2mu}\sigma\big)^2\mspace{-2mu}a d\mspace{-1.5mu}+\mspace{-2mu}\frac{(K\mspace{-2mu}+\mspace{-2mu}\sigma)^4 a^2d}{1\mspace{-2mu}-\mspace{-2mu}2(K\mspace{-2mu}+\mspace{-2mu}\sigma)^2 a}\right)\mspace{-3mu},
\end{align*}
Setting $a=\frac{1}{4(K+\sigma)^2}$ and combining \eqref{EQ:W1_CS_bound}-\eqref{EQ:W1_int2_bound} yields
\begin{equation}
\EE_{\prodmu}\Gwass(\empmu,\mu)\leq c_1^d \sigma\sqrt{2d}\left(1+\frac{K}{\sigma}\right)^{\frac{d}{2}+1}e^{\frac{3d}{16}}\frac{1}{\sqrt{n}},\label{EQ:W1_bound_constant}
\end{equation}
where $c_1$ is the constant from Lemma \ref{LEMMA:density_bound}. We note that a better constant can be achieved by assuming $\Noise=\Gauss$ \cite{Goldfeld2019convergence}, but we chose to sacrifice that in favor of generality.


\subsection{Proof of Corollary \ref{COR:gwass_emp_whp}}\label{SUBSEC:gwass_emp_whp_proof}

The main tool used in this proof is McDiarmid's inequality:
\vspace{2mm}
\begin{lemma}[McDiarmid's Inequality]
Let $X^n\triangleq(X_1,\ldots,X_n)$ be an $n$-tuple of $\cX$-valued independent random variables. Suppose $g:\cX^n\to\RR$ is a map that for any $i=1,\ldots,n$ and $x_1,\ldots,x_n,x_i'\in\cX$ satisfies
\vspace{2mm}
\begin{equation}
    \big|g(x^n)-g(x_1,\ldots,x_{i-1},x'_i,x_{i+1},\ldots,x_n)\big|\leq c_i,\label{EQ:bdd_diff}
\end{equation}
for some non-negative $\{c_i\}_{i=1}^n$. Then for any $t>0$:
\begin{subequations}
\begin{align}
    \PP\Big(g(X^n)-\EE g(X^n)\geq t\Big)&\leq e^{-\frac{2t^2}{\sum_{i=1}^nc_i^2}}\label{EQ:McDiarmid1}\\
    \PP\Big(\big|g(X^n)-\EE g(X^n)\big|\geq t\Big)&\leq 2e^{-\frac{2t^2}{\sum_{i=1}^nc_i^2}}\label{EQ:McDiarmid2}
\end{align}
\end{subequations}
\end{lemma}

Let $g(X^n)\triangleq\Gwass(\empmu,\mu)$ and use Kantorovich-Rubinstein duality:
\begin{align*}
    g(X^n)&=\sup_{\|f\|_{\mathsf{Lip}}\leq 1} \EE_{\empmu\ast\Noise}f-\EE_{\mu\ast\Noise}f\\
    &=\sup_{\|f\|_{\mathsf{Lip}}\leq 1} \frac{1}{n}\sum_{i=1}^n(f\ast\noise)(X_i)-\EE_{\mu}\big[f\ast\noise\big].
\end{align*}

Fix $i\in\{1,\ldots,n\}$ and $x_1,\ldots,x_n,x_i'\in\cX$. Property \eqref{EQ:bdd_diff} follows by first observing that:
\begin{align*}
    &n\Big(g(x^n)-g(x_1,\ldots,x_{i-1},x'_i,x_{i+1},\ldots,x_n)\Big)\\
    &=\mspace{-3mu}\sup_{\|f\|_{\mathsf{Lip}}\leq 1}\mspace{-3mu}\Bigg\{\mspace{-1mu}\mspace{-2mu}\sum_{j\neq i}(f\mspace{-2mu}\ast\mspace{-2mu}\noise)(x_j)\mspace{-2mu}-\mspace{-2mu}\EE_{\mu}\big[\mspace{-1mu}f\mspace{-2mu}\ast\mspace{-2mu}\noise\mspace{-1mu}\big]\mspace{-3mu}+\mspace{-3mu}(f\mspace{-2mu}\ast\mspace{-2mu}\noise)(x_i)\Bigg\}\\
    &-\sup_{\|h\|_{\mathsf{Lip}}\leq 1}\mspace{-3mu}\Bigg\{\mspace{-1mu}\mspace{-2mu}\sum_{j\neq i}(h\mspace{-2mu}\ast\mspace{-2mu}\noise)(x_j)\mspace{-2mu}-\mspace{-2mu}\EE_{\mu}\big[\mspace{-1mu}h\mspace{-2mu}\ast\mspace{-2mu}\noise\mspace{-1mu}\big]\mspace{-3mu}+\mspace{-3mu}(h\mspace{-2mu}\ast\mspace{-2mu}\noise)(x_i')\Bigg\}\\
    &\leq \sup_{\|f\|_{\mathsf{Lip}}\leq 1}(f\ast\noise)(x_i)-(f\ast\noise)(x_i').\numberthis\label{EQ:McDiarmid_UB}
\end{align*}
Then we note that Lipschitzness of $f$ implies that $f\ast\noise$ is also Lipschitz.

\begin{lemma}[Lipschitz after Convolution]\label{LEMMA:Lip_conv_Lip}
If $f:\RR^d\to\RR$ has $\|f\|_{\mathsf{Lip}}\leq L$, then $\|f\ast g\|_{\mathsf{Lip}}\leq L$ for any PDF $g:\RR^d\to\RR_{\geq 0}$.
\end{lemma}
The proof is immediate and thus omitted. Combining Lemma \ref{LEMMA:Lip_conv_Lip} and \eqref{EQ:McDiarmid_UB}, we obtain
\begin{equation*}
    \Big|g(x^n)-g(x_1,\ldots,x_{i-1},x'_i,x_{i+1},\ldots,x_n)\Big|\leq \frac{\mathsf{diam}(\cX)^2}{n},
\end{equation*}
for all $i=1,\ldots,n$ and $x_1,\ldots,x_n,x_i'\in\cX$.

Applying McDiarmiad's inequality \eqref{EQ:McDiarmid2} for $g(X^n)=\Gwass(\empmu,\mu)$ produces \eqref{EQ:gwass_emp_whp_equation}. Taking $t=\Theta\left(\frac{\log n}{\sqrt{n}}\right)$ and inserting into \eqref{EQ:McDiarmid1} gives \eqref{EQ:gwass_emp_whp_equation1}.


\section{Summary and Concluding Remarks}\label{SEC:summary}

We proposed a novel Gaussian-smoothed framework for OT defined as $\gwass(\mu,\nu)\triangleq\wass(\mu\ast\Gauss,\nu\ast\Gauss$). This GOT distance was shown to inherit the metric structure (and the metrization of weak convergence) from the regular 1-Wasserstein distance. As a function of $\sigma$, $\gwass(\mu,\nu)$ is a continuous and monotonically decreasing function maximized at $\wass^{(0)}(\mu,\nu)=\wass(\mu,\nu)$. Furthermore, as $\gwass(\mu,\nu)\xrightarrow[\sigma\to 0]{}\wass(\mu,\nu)$, optimal transport plans for $\gwass(\mu,\nu)$ weakly converge to an optimal plan for $\wass(\mu,\nu)$. Finally, we explored statistical properties of $\gwass$, studying the convergence rate of $\EE\gwass(\empmu,\mu)$ to 0, where $\empmu$ is the empirical measure induced by $n$ i.i.d. samples from $\mu$. Building on \cite{Goldfeld2019convergence}, we showed that $\wass(\empmu\ast\Noise,\mu\ast\Noise)\in O\left(n^{-\frac{1}{2}}\right)$ in all dimensions, for any subgaussian noise distribution $\Noise$ with a monotone and bounded density. In particular, $\gwass$ alleviates the curse of dimensionality in the one-sample (and hence also in the weaker two-sample) regime. This stands in striking contrast to classic 1-Wasserstein distance, which converge at most as $n^{-\frac{1}{d}}$, while no results are available for entropic OT with distance cost. 
These theoretical findings were verified through an empirical study, posing GOT as an appealing alternative to the popular entropically regularized OT methods.


Attractive next steps include the design of efficient algorithms tailored for GOT computation. While any method for computing Wasserstein distances is also applicable for GOT, it possesses additional structure one may exploit. We plan to leverage this structure in our future algorithmic designs, and explore avenues for their use in generative modeling and other OT applications. Additional directions include examining alternative noise models and their comparison to the Gaussian-smoothed framework.

\bibliographystyle{IEEEtran}
\bibliography{ref}

\begin{thebibliography}{10}
\providecommand{\url}[1]{#1}
\csname url@samestyle\endcsname
\providecommand{\newblock}{\relax}
\providecommand{\bibinfo}[2]{#2}
\providecommand{\BIBentrySTDinterwordspacing}{\spaceskip=0pt\relax}
\providecommand{\BIBentryALTinterwordstretchfactor}{4}
\providecommand{\BIBentryALTinterwordspacing}{\spaceskip=\fontdimen2\font plus
\BIBentryALTinterwordstretchfactor\fontdimen3\font minus
  \fontdimen4\font\relax}
\providecommand{\BIBforeignlanguage}[2]{{%
\expandafter\ifx\csname l@#1\endcsname\relax
\typeout{** WARNING: IEEEtran.bst: No hyphenation pattern has been}%
\typeout{** loaded for the language `#1'. Using the pattern for}%
\typeout{** the default language instead.}%
\else
\language=\csname l@#1\endcsname
\fi
#2}}
\providecommand{\BIBdecl}{\relax}
\BIBdecl

\bibitem{kantorovich1942translocation}
L.~V. Kantorovich, ``On the translocation of masses,'' in \emph{USSR Academy of
  Science (Doklady Akademii Nauk USSR}, vol.~37, 1942, pp. 199--201.

\bibitem{ho2017multilevel}
N.~Ho, X.~L. Nguyen, M.~Yurochkin, H.~H. Bui, V.~Huynh, and D.~Phung,
  ``Multilevel clustering via wasserstein means,'' in \emph{International
  Conference on Machine Learning (ICML-2017)}, Sydney, Australia, Jul. 2017,
  pp. 1501--1509.

\bibitem{iyer2014maximum}
A.~Iyer, S.~Nath, and S.~Sarawagi, ``Maximum mean discrepancy for class ratio
  estimation: Convergence bounds and kernel selection,'' in \emph{International
  Conference on Machine Learning}, 2014, pp. 530--538.

\bibitem{courty2016optimal}
N.~Courty, R.~Flamary, D., Tuia, and A.~Rakotomamonjy, ``Optimal transport for
  domain adaptation,'' \emph{IEEE Transactions on Pattern Analysis and Machine
  Intelligence}, vol.~39, no.~9, pp. 1853--1865, Oct. 2016.

\bibitem{courty2014domain}
N.~Courty, R.~Flamary, and D.~Tuia, ``Domain adaptation with regularized
  optimal transport,'' in \emph{European Conference on Machine Learning and
  Knowledge Discovery in Databases (ECML PKDD 2014)}, Nancy, France, Sep. 2014,
  pp. 274--289.

\bibitem{arjovsky2017wasserstein}
M.~Arjovsky, S.~Chintala, and L.~Bottou, ``Wasserstein generative adversarial
  networks,'' in \emph{International Conference on Machine Learning
  (ICML-2017)}, Sydney, Australia, Jul. 2017, pp. 214--223.

\bibitem{gulrajani2017improved}
I.~Gulrajani, F.~Ahmed, M.~Arjovsky, V.~Dumoulin, and A.~C. Courville,
  ``Improved training of {W}asserstein {GAN}s,'' in \emph{Advances in Neural
  Information Processing Systems (NeurIPS-2017)}, Long Beach, CA, US, Dec.
  2017, pp. 5767--5777.

\bibitem{rubner2000earth}
Y.~Rubner, C.~Tomasi, and L.~J. Guibas, ``The earth mover's distance as a
  metric for image retrieval,'' \emph{International Journal of Computer
  Vision}, vol.~40, no.~2, pp. 99--121, Nov. 2000.

\bibitem{sandler2011nonnegative}
R.~Sandler and M.~Lindenbaum, ``Nonnegative matrix factorization with earth
  mover's distance metric for image analysis,'' \emph{IEEE Transactions on
  Pattern Analysis and Machine Intelligence}, vol.~33, no.~8, pp. 1590--1602,
  Jan. 2011.

\bibitem{li2013novel}
P.~Li, Q.~Wang, and L.~Zhang, ``A novel earth mover's distance methodology for
  image matching with gaussian mixture models,'' in \emph{IEEE International
  Conference on Computer Vision (ICCV-2013)}, Sydney, Australia, Dec. 2013, pp.
  1689--1696.

\bibitem{alvarez2018gromov}
D.~Alvarez-Melis and T.~S. Jaakkola, ``{Gromov-Wasserstein} alignment of word
  embedding spaces,'' \emph{arXiv preprint arXiv:1809.00013}, Aug. 2018.

\bibitem{yurochkin2019hierarchical}
M.~Yurochkin, S.~Claici, E.~Chien, F.~Mirzazadeh, and J.~Solomon,
  ``Hierarchical optimal transport for document representation,'' \emph{arXiv
  preprint arXiv:1906.10827}, Jun 2019.

\bibitem{grave2019unsupervised}
E.~Grave, A.~Joulin, and Q.~Berthet, ``Unsupervised alignment of embeddings
  with {Wasserstein} procrustes,'' in \emph{International Conference on
  Artificial Intelligence and Statistics (AISTATS-2019)}, Okinawa, Japan, Apr.
  2019, pp. 1880--1890.

\bibitem{panaretos2019statistical}
V.~M. Panaretos and Y.~Zemel, ``Statistical aspects of wasserstein distances,''
  \emph{Annual Review of Statistics and its Application}, vol.~6, pp. 405--431,
  Mar. 2019.

\bibitem{villani2008optimal}
C.~Villani, \emph{Optimal transport: old and new}.\hskip 1em plus 0.5em minus
  0.4em\relax Springer Science \& Business Media, 2008, vol. 338.

\bibitem{varadarajan1958convergence}
V.~S. Varadarajan, ``On the convergence of sample probability distributions,''
  \emph{Sankhy{\=a}: The Indian Journal of Statistics (1933-1960)}, vol.~19,
  no. 1/2, pp. 23--26, Feb 1958.

\bibitem{dudley1969speed}
R.~M. Dudley, ``The speed of mean {Glivenko-Cantelli} convergence,'' \emph{Ann.
  Math. Stats.}, vol.~40, no.~1, pp. 40--50, Feb. 1969.

\bibitem{dobric1995asymptotics}
V.~Dobri{\'c} and J.~E. Yukich, ``Asymptotics for transportation cost in high
  dimensions,'' \emph{J. Theoretical Prob.}, vol.~8, no.~1, pp. 97--118, Jan.
  1995.

\bibitem{cuturi2013sinkhorn}
M.~Cuturi, ``Sinkhorn distances: {Lightspeed} computation of optimal
  transport,'' in \emph{Advances in Neural Information Processing Systems
  (NeurIPS-2013)}, Stateline, NV, US, Dec. 2013, pp. 2292--2300.

\bibitem{altschuler2017near}
J.~Altschuler, J.~Weed, and P.~Rigollet, ``Near-linear time approximation
  algorithms for optimal transport via {Sinkhorn} iteration,'' in
  \emph{Advances in Neural Information Processing Systems (NeurIPS-2017)}, Long
  Beach, CA, US, Dec. 2017, pp. 1964--1974.

\bibitem{genevay2016stochastic}
A.~G. M., Cuturi, G.~Peyr{\'e}, and F.~Bach, ``Stochastic optimization for
  large-scale optimal transport,'' in \emph{Advances in Neural Information
  Processing Systems (NeurIPS-2016)}, Barcelona, Spain, Dec. 2017, pp.
  3440--3448.

\bibitem{montavon2016wasserstein}
G.~Montavon, K.-R. M{\"u}ller, and M.~Cuturi, ``Wasserstein training of
  restricted boltzmann machines,'' in \emph{Advances in Neural Information
  Processing Systems (NeurIPS-2016)}, Barcelona, Spain, Dec. 2016, pp.
  3718--3726.

\bibitem{rigollet2018entropic}
P.~Rigollet and J.~Weed, ``Entropic optimal transport is maximum-likelihood
  deconvolution,'' \emph{Comptes Rendus Mathematique}, vol. 356, no. 11-12, pp.
  1228--1235, Nov 2018.

\bibitem{genevay2019sample}
A.~Genevay, L.~Chizat, F.~Bach, M.~Cuturi, and G.~Peyr{\'e}, ``Sample
  complexity of {Sinkhorn} divergences,'' in \emph{International Conference on
  Artificial Intelligence and Statistics (AISTATS-2019)}, Okinawa, Japan, Apr.
  2019, pp. 1574--1583.

\bibitem{mena2019statistical}
G.~Mena and J.~Weed, ``Statistical bounds for entropic optimal transport:
  sample complexity and the central limit theorem,'' \emph{arXiv preprint
  arXiv:1905.11882}, May 2019.

\bibitem{feydy2018interpolating}
J.~Feydy, T.~S{\'e}journ{\'e}, F.-X. Vialard, S.-I. Amari, A.~Trouv{\'e}, and
  G.~Peyr{\'e}, ``Interpolating between optimal transport and mmd using
  sinkhorn divergences,'' \emph{arXiv preprint arXiv:1810.08278}, Oct. 2018.

\bibitem{bigot2019central}
J.~Bigot, E.~Cazelles, and N.~Papadakis, ``Central limit theorems for
  entropy-regularized optimal transport on finite spaces and statistical
  applications,'' \emph{arXiv preprint arXiv:1711.08947}, 2019.

\bibitem{dalmaso2012Gamma_conv}
G.~D. Maso, \emph{An introduction to $\Gamma$-convergence}.\hskip 1em plus
  0.5em minus 0.4em\relax Springer Science \& Business Media, 2012, vol.~8.

\bibitem{Goldfeld2019convergence}
Z.~Goldfeld, K.~Greenewald, Y.~Polyanskiy, and J.~Weed, ``Convergence of
  smoothed empirical measures with applications to entropy estimation,''
  \emph{arXiv preprint arXiv:1905.13576}, May 2019.

\bibitem{hsu2012tail}
D.~Hsu, S.~Kakade, and T.~Zhang, ``A tail inequality for quadratic forms of
  subgaussian random vectors,'' \emph{Electronic Communications in
  Probability}, vol.~17, 2012.

\bibitem{vershynin2010introduction}
R.~Vershynin, ``Introduction to the non-asymptotic analysis of random
  matrices,'' \emph{arXiv preprint arXiv:1011.3027}, 2010.

\end{thebibliography}
	
\clearpage

\newpage

\appendix
\section*{Appendix}
\section{Proof of Lemma \ref{LEMMA:density_bound}}\label{APPEN:density_bound_proof}

Recall that $\noise(t) = \prod_{j=1}^d \noised(t_j)$, where $\noised$ is $\sigma$-subgaussian, zero mean, bounded, and monotonically decreasing as $t_j$ moves away from zero. We first analyze the one-dimensional densities $\noised$, and show that there exists a constant $c>0$, such that
\begin{equation}
    \noised(t) \leq c e^{2 \delta |t| - \delta^2 -\log \delta} \gaussd(t),\quad \forall t\in\RR,\label{EQ:1d_bound}
\end{equation}
where $\gaussd$ is a scalar Gaussian density (zero mean and $\sigma^2$ variance). We prove \eqref{EQ:1d_bound} for $t>0$; the $t < 0$ case is identical. 

Note that the $\sigma$-subgaussianity of $\noised$ (Def. \ref{DEF:SG}) implies that
\begin{equation}
\mathbb{E}_{\noised}\left[e^{\alpha X}\right] \leq e^{\frac{1}{2} \sigma^2 \alpha^2},\quad\forall \alpha\in\RR,\label{EQ:subGNoise}
\end{equation}
which by \cite{vershynin2010introduction} yields
\begin{equation}
    \PP_{\noised}\big((-\infty,t)\cup(t,\infty)\big)\leq \exp(1-t^2/(2\sigma^2)) = c'\gaussd(t),
\end{equation}
where $c'=\sqrt{2\pi\sigma^2e^2}$. Consequently, for any $t^\star$, 
\begin{align*}
\mathbb{P}_{\noised}\big((t^\star-\delta,t^\star]\big) &\leq \mathbb{P}_{\noised}\big((t^\star-\delta,\infty)\big)\\
&\leq c' \gaussd(t^\star-\delta)\\
&= c' e^{(t^\star)^2 - (t^\star-\delta)^2}\gaussd(t^\star)\\
&= c' e^{2\delta t^\star - \delta^2} \gaussd(t^\star).\numberthis \label{eq:sgbd}
\end{align*}
Now, since $\tilde{g}_\sigma(t)$ monotonically decreases as $t$ moves away from zero, for any $t^\star \geq \delta$ we have $\mathbb{P}_{\noised}\big((t^\star-\delta,t^\star]\big) \geq \delta \tilde{g}_\sigma(t^\star)$. Substituting this into \eqref{eq:sgbd}, we have for all $t^\star \geq \delta$ that
\begin{align*}
\delta \tilde{g}_\sigma(t^\star) &\leq c' e^{2\delta t^\star - \delta^2} \gaussd(t^\star),\\
\tilde{g}_\sigma(t^\star) &\leq c' e^{2\delta t^\star - \delta^2-\log \delta} \gaussd(t^\star).
\end{align*}
Repeating the argument for $t < 0$ then yields 
\[
\tilde{g}_\sigma(t) \leq c' e^{2\delta |t| - \delta^2-\log \delta} \gaussd(t)
\]
for all $|t| \geq \delta$.
Since $\tilde{g}_\sigma$ is bounded, $\sup_{|t|\leq \delta}\noised(t)\left(e^{2 \delta t - \delta^2-\log \delta} \gaussd(t)\right)^{-1}$ exists, and hence \eqref{EQ:1d_bound} holds (for all $t \in \mathbb{R}$) with 
\[
c = \max\left[c',\sup_{|t|\leq \delta}\noised(t)\left(e^{2 \delta t - \delta^2-\log \delta} \gaussd(t)\right)^{-1}\right].
\]

Extending to the full $d$-dimensional distribution, note that since $t^2 + 1 > |t|$ for all $t$, we have that $\noised(t) \leq c e^{2 \delta t^2 + 2\delta - \delta^2 -\log \delta}\gaussd(t)$ for all $t$. We can then write
\begin{equation}
    \noise(t) \leq (c')^d e^{2 \delta \|t\|^2 + 2d\delta - d\delta^2 -d\log \delta} \gauss(t),\label{EQ:lemma_final}
\end{equation}
which establishes the lemma after collecting terms.



\end{document}